\documentclass[11pt]{amsart}

\usepackage[T1]{fontenc}
\usepackage[active]{srcltx}
\usepackage[leqno]{amsmath}
\usepackage{amssymb}
\usepackage{amsfonts}
\usepackage{amsthm}
\usepackage{latexsym}
\usepackage{euscript}
\usepackage{mathdots}
\usepackage{oldgerm}
\usepackage{cite}
\usepackage{mathtools}

\usepackage{pifont}
\usepackage{mdwlist}
\usepackage{mathrsfs}
\usepackage{comment}
\usepackage{setspace}
\usepackage{mathdots}
\usepackage{tikz}
\usepackage{tikz-cd}
\usepackage{graphicx}%,pst-all,bm}
\usepackage{fixltx2e}
\usepackage{flushend}
\usepackage{subfig}
\usepackage{enumerate}
\usepackage{caption}
\usepackage{longtable}
\usepackage{verbatim}
\usepackage{relsize}
\usepackage{exscale}
\usepackage{hyperref}

\usepackage{graphicx}
\usepackage{tikz}
\usetikzlibrary{patterns,angles,quotes,arrows,shapes,calc,positioning}
\usepackage[europeanresistors,europeaninductors,siunitx,nooldvoltagedirection]{circuitikz}

\numberwithin{equation}{section}

\newcommand{\cC}{\mathcal{C}}

\newcommand{\cP}{\mathcal{P}}

\newcommand{\cS}{\mathcal{S}}

\newcommand{\cX}{\mathcal{X}}
\newcommand{\cY}{\mathcal{Y}}

\makeindex             % used for the subject index
\usetikzlibrary{arrows}

\theoremstyle{theorem}
\newtheorem{Prop}{Proposition}[section]

\newtheorem{Thm}[Prop]{Theorem}

\theoremstyle{definition}
\newtheorem{Def}[Prop]{Definition}
\newtheorem{Rem}[Prop]{Remark}
\newtheorem{Rems}[Prop]{Remarks}

\newtheorem{Ex}[Prop]{Example}

\newcommand{\tto}{\!\rightarrow\!}
\newcommand{\N}{{\mathbb{N}}}

\newcommand{\R}{{\mathbb{R}}}
\DeclareMathOperator{\rk}{rk}

\DeclareMathOperator{\im}{im}

\newcommand{\del}{\partial}
\newcommand{\setdef}[2]{\left\{\ #1\ \left|\ \vphantom{#1} #2\ \right.\right\}}

\newcommand{\ds}[1]{{\rm \, d} #1 \,}

{\left(\begin{smallmatrix}}%            begin code
{\end{smallmatrix}\right)}%             end code

{\left[\begin{smallmatrix}}%            begin code
{\end{smallmatrix}\right]}%             end code

\renewcommand{\restriction}{\mathord{\upharpoonright}}

\begin{document}
\title[On a global implicit function theorem]{On an extension of a global implicit function theorem}

\author[T.~Berger]{Thomas Berger}
\address{Institut f\"ur Mathematik, Universit\"at Paderborn, Warburger Str.~100, 33098~Paderborn, Germany}
\email{thomas.berger@math.upb.de}

\author[F.~E.~Haller]{Fr\'ed\'eric E. Haller}
\address{Fachbereich Mathematik, Universit\"{a}t Hamburg, Bundesstrasse 55, D-20146 Hamburg, Germany}
\email{frederic.haller@uni-hamburg.de}
%\thanks{}

\date{\today}


\begin{abstract}
We study the existence of global implicit functions for equations defined on open subsets of Banach spaces. The partial derivative with respect to the second variable is only required to have a left inverse instead of being invertible. Generalizing known results, we provide sufficient criteria which are easy to check. These conditions essentially rely on the existence of diffeomorphisms between the respective projections of the set of zeros and appropriate Banach spaces, as well as a corresponding growth bound. The projections further allow to consider cases where the global implicit function is not defined on all of the open subset corresponding to the first variable.
\end{abstract}

\maketitle

%\noindent{
%{\bf Nomenclature:}\\[2mm]
%\begin{tabular}[ht]{p{43pt}p{187pt}} %265
%$\N$, $\N_0$
%%, $\Z$
%&  set of natural numbers, $\N_0=\N\cup\{0\}$\\[1mm]
%$\C_{+}$
%   & open set of complex numbers with positive real part \\[1mm]
%$\R^{n\times m}$
%  &   the set of real $n\times m$ matrices\\[1mm]
%$\rk A,\, \im A$
%  &   rank and image of $A\in \R^{n\times m}$\\[1mm]
%$\Gl_n(\R)$ & the group of invertible matrices in $\R^{n\times n}$\\[1mm]
%$M >_{\cV} 0$ & $:\Longleftrightarrow\  \forall\, x\in\cV\setminus\{0\}:\ x^\top M x > 0$, for a matrix $M\in\R^{n\times n}$ and a subspace $\cV\subseteq\R^n$; similarly, the notation $M =_{\cV} 0$ is introduced\\[1mm]
%%$I_n$ && identity matrix of size $n\times n$\\[2mm]
%%$A^\dagger$ &$=(A^\top A)^{-1} A^\top$ for $A\in\R^{n\times m}$ with $\rk A=m$\\[2mm]
%$\cC^k(X\tto Y)$ & set of $k$-times continuously differentiable functions $f:X\tto Y$, $k\in\N_0\cup\{\infty\}$; $\cC(X\tto Y)\coloneqq\cC^0(X\tto Y)$; if $k=\infty$ the function $f$ is called \emph{smooth}\\[1mm]
%$\dom\, f$ & the domain of the function $f$\\[1mm]
%$\left.f\right|_I$ & restriction of the function $f$ to the set $I$
%\end{tabular}
%}

%%%%%%%%%%%%%%%%%%%%%%%%%%%%%%%%%%%%%%%%%%%%%%%%%%%%%%%%%%%%%%%%%%%%%%%%%%%%%%%%%%%%%%%%%%%%%%%%%%%%%%%%%

%\newpage
%
%%%%%%%%%%%%%%%%%%%%%%%%%%%%%%%%%%%%%%%%%%%%%%%%%%%%%%%%%%%%%%%%%%%%%%%%%%%%%%%%%%%%%%%%%%%%%%%%%%%%%%%%%%%%%
\section{Introduction}\label{Sec:Intr}
%%%%%%%%%%%%%%%%%%%%%%%%%%%%%%%%%%%%%%%%%%%%%%%%%%%%%%%%%%%%%%%%%%%%%%%%%%%%%%%%%%%%%%%%%%%%%%%%%%%%%%%%%%%%%
%

When dealing with nonlinear dynamical systems with constraints, i.e., implicit differential equations of the form
\[
    F\big( x(t), \dot x(t)\big) = 0,
\]
where the number of equations does not match the number of variables, it is often necessary to solve this equation for~$\dot x(t)$, preferably globally in the form $\dot x(t) = g\big(x(t)\big)$. This means to find a global implicit function of the equation $F(x,y) = 0$. Numerous results on global implicit function theorems exist, and we mention the relevant literature. However, most results involve conditions which are not easy to check in practice. In the present paper, we provide a novel extension of the global implicit function theorem under conditions which can easily be verified.

In the following we summarize some results on global implicit functions, tailored to be applicable in our framework. We consider equations of the form $F(x,y)=0$ for which we want to find a unique maximal solution $y(x)$. There are several approaches available in the literature which provide a solution to this problem, see e.g.~\cite{Rhei69} for an early result. Most works concentrate on the case that the partial derivative $\frac{\del F}{\del y}(x,y)$ is invertible for all~$(x,y)$, i.e., $F(x,y)=0$ is locally solvable for $y(x)$ in a neighborhood of every point~$(a,b)$ such that $F(a,b)=0$. We discuss some important work:
\begin{itemize}
  \item For $F:X\times Y\to \R^l$, where $X\subseteq\R^m$ and $Y\subseteq\R^n$ are open and $X$ is convex, Sandberg~\cite{Sand81} provides necessary and sufficient conditions for the existence of a unique $g\in\cC(X, Y)$ such that $F^{-1}(0) = \setdef{(x,g(x))}{x\in X}$. However, the conditions are not easy to check; in particular, it needs to be guaranteed that
      \begin{equation}\label{eq:Sand}
      \begin{aligned}
        & \text{for some $x_0\in X$ there exists exactly one}\\
        & \text{$y_0\in Y$ such that $F(x_0,y_0)=0$.}
      \end{aligned}
       \end{equation}
       Furthermore, the result of Sandberg is not applicable in the case that the maximal solution~$g$ is not defined on all of~$X$.
  \item Using the theory of covering maps, Ichiraku~\cite{Ichi85} improves the characterization of Sandberg. Nevertheless, the condition~\eqref{eq:Sand} is still present and the results are only applicable in the case of globally defined~$g$. However, in~\cite[Thm.~5]{Ichi85} it is shown that in the case $X=\R^m$, $Y=\R^n$ and $l=m$ for the existence of a unique solution $g\in\cC(\R^m,\R^n)$ it is sufficient that $\frac{\del F}{\del y}(x,y)$ is invertible for all~$(x,y)\in\R^m\times\R^n$, condition~\eqref{eq:Sand} holds and
      \begin{equation}\label{eq:Ichi}
        \forall\, (x,y)\in F^{-1}(0):\ \ \left\|\left(\frac{\del F}{\del y}(x,y)\right)^{-1}\right\|\cdot \left\| \frac{\del F}{\del x} (x,y) \right\| \le M
      \end{equation}
      for some $M\ge 0$.
  \item The above result of Ichiraku has in turn be improved by Gut\'{u} and Jaramillo~\cite[Cor.~5.3]{GutuJara07}, who showed that the condition~\eqref{eq:Sand} can be replaced by the intuitive condition ``$F^{-1}(0)$ is connected'' and in the condition~\eqref{eq:Ichi} the constant~$M$ can be replaced by the term $\omega(\|y\|)$, where $\omega:[0,\infty)\to(0,\infty)$ is a continuous \emph{weight}, which means that $\omega$ is nondecreasing and
      \[
        \int_0^\infty \frac{\ds{t}}{\omega(t)} = \infty.
      \]
      These conditions are indeed easy to check. The only drawback is that~$F$ needs to be defined on all of $\R^m\times \R^n$ and the solution~$g$ is defined on all of~$\R^m$.
  \item A result which is similar to that of Gut\'{u} and Jaramillo, but holds for some $X\subseteq\R^m$ which is open, connected and starlike with respect to some $a\in X$ such that $F(a,b)=0$ for some $b\in Y=\R^n$, has been derived by Cristea~\cite{Cris07}. The assumption of connectedness of $F^{-1}(0)$ is not needed, however a version of assumption~\eqref{eq:Ichi} (with $M=\omega(\|y\|)$) is required to hold on all of $X\times \R^n$.
  \item Yet another approach has been pursued by Zhang and Ge~\cite{ZhanGe06} who show that for existence of a unique solution $g\in\cC(\R^m,\R^n)$ it is sufficient that the element-wise absolute value of $\frac{\del F}{\del y}$ is uniformly strictly diagonally dominant in the sense that there exists $d>0$ such that
      \[
        \left| \left(\frac{\del F}{\del y}(x,y)\right)_{ii}\right| - \sum_{j\ne i} \left| \left(\frac{\del F}{\del y}(x,y)\right)_{ij}\right| \geq d
      \]
      for all $(x,y)\in\R^m\times \R^n$ and all $i=1,\ldots,n$. While this condition is easy to check, it is very restrictive as it already excludes a lot of linear equations $Ax + By = 0$ where~$B$ is not strictly diagonally dominant, but invertible.
\end{itemize}

As discussed above, typical limitations of the approaches are that~$F$ needs to be defined on all of $\R^m\times \R^n$ or the solution~$g$ is required to be globally defined. In~\cite{Blot91} these limitations are resolved as~$X$ and~$Y$ are assumed to be open and~$X$ is connected, and maximal solutions of $F(x,y)=0$ are considered in every connected component of $F^{-1}(0)$. Assuming that $Z\coloneqq F^{-1}(0)$ is connected we may then find a solution $g\in\cC(\pi_1(Z), Y)$, where $\pi_1:X\times Y\to X,\ (x,y)\mapsto x$ is the projection onto the first component, provided that $\pi_1(Z)$ is open and simply connected and $\pi_1:Z\tto \pi_1(Z)$ ``lifts lines'' (for a precise definition see~\cite[Def.~1.1]{Plas74}). This result can be extended in a straightforward way to the case where $l\ge m$ and $\rk \frac{\del F}{\del y} (x,y) = n$ for all $(x,y)\in X\times Y$ since it is only necessary to show that~$\pi_1$ is locally a homeomorphism, which replaces the condition that $F(x,y)=0$ is locally solvable for~$y(x)$ as in~\cite[Thm.~4]{Blot91}; then~\cite[Lem.~1]{Blot91} can still be applied to $\pi_1:Z\tto \pi_1(Z)$. The drawback of this result is that the condition ``$\pi_1:Z\tto \pi_1(Z)$ lifts lines'' is not easy to check.

In the present paper, we provide a generalization of~\cite[Cor.~5.3]{GutuJara07} to the case of functions defined only on open subsets and where the partial derivative $\frac{\del F}{\del y}$ is only required to have a left inverse instead of being invertible. The crucial assumption is that the projections $\pi_i(Z)$ on the $i$th component, $i=1,2$, are diffeomorphic to some Banach spaces and the transformation of the equation $F(x,y)=0$ satisfies a generalized version of~\eqref{eq:Ichi}. We stress that this assumption in particular implies that~$\pi_i(Z)$ must be open and simply connected. The main result is presented in Section~\ref{Sec:main} and a discussion together with some illustrative examples is given in Section~\ref{Sec:Ex}.

%%%%
\section{Main result}\label{Sec:main}
%%%

In this section we state and prove the following main result of the paper.

\begin{Thm}\label{Thm:global-IF}
  Let $X\subseteq \mathcal{U},\ Y\subseteq \mathcal{V}$ be open sets, $\mathcal{U,V,Z}$ be Banach spaces, $F\in\cC^1(X\times Y,\mathcal{Z})$ and
  \[
     Z \subseteq \setdef{ (x,y)\in X\times Y }{ F(x,y) = 0}
  \]
  be such that
  \begin{enumerate}
    \item $Z$ is path-connected and closed in $X\times Y$;
    \item $\forall\, (x,y)\in Z\  \exists\, S(x,y)\in\mathcal{L}(\mathcal{Z},\mathcal{V}):\ S(x,y)D_yF(x,y)={\rm id}_\mathcal{V}$;\footnote{Here $\mathcal{L}(\mathcal{Z},\mathcal{V})$ denotes the Banach space of all bounded linear operators $A:\mathcal{Z}\to\mathcal{V}$ and ${\rm id}_\mathcal{V}:\mathcal{V}\to\mathcal{V},\ v\mapsto v$ is the identity operator on~$\mathcal{V}$.}
    \item for the projections $\pi_i(p_1,p_2)=p_i$ with $i\in\{1,2\}$, $(p_1,p_2)\in\mathcal{U}\times\mathcal{V}$, there exist diffeomorphisms $\phi:\pi_1(Z)\tto\mathcal{X}$, $\psi:\pi_2(Z)\tto\mathcal{Y}$ for some Banach spaces $\mathcal{X},\mathcal{Y}$, and a continuous weight $\omega:[0,\infty)\to(0,\infty)$ such that for all $(x,y)\in Z$ we have
    % \[
    %     T(x,y)\coloneqq \left( \left(\frac{\del F}{\del y} (x,y)\right)^\top \frac{\del F}{\del y} (x,y)\right)^{-1} \left(\frac{\del F}{\del y} (x,y)\right)^\top,
    % \]
    % that

\[
	\left\| D\psi({y})\cdot{S}({x},{y})\right\|_{\mathcal{L}(\mathcal{Z},\mathcal{Y})} \cdot\left\|D_{{x}}{F}({x},{y})\cdot \big(D\phi({x})\big)^{-1}\right\|_{\mathcal{L}(\mathcal{X},\mathcal{Z})} \le \omega(\|\psi(y)\|_\mathcal{Y}).
\]

\end{enumerate}
Then there exists a unique $g\in\cC(\pi_1(Z), Y)$ such that
\[
    \setdef{(x,g(x))}{x\in \pi_1(Z)} = Z,
\]
and $g$ is Fr\'echet-differentiable at every $x\in\pi_1(Z)$.
\end{Thm}

The proof of Theorem~\ref{Thm:global-IF} requires us to recall the following concepts, which can be found in \cite[pp.~77--80]{GutuJara07}.

\begin{Def}
Let $Z$ be a metric space, and let $\mathcal{P}$ be a family of continuous paths in $Z$. We say that $Z$ is \emph{$\mathcal{P}$-connected}, if the following conditions hold:
\begin{enumerate}[(1)]
\item If the path $p:[a,b]\to Z$ belongs to $\mathcal{P}$, then the reverse path $\overline{p}$, defined by $\overline{p}(t)=p(a-t+b)$, also belongs to $\mathcal{P}$;
\item Every two points in $Z$ can be joined by a path in $\mathcal{P}$.
\end{enumerate}
We say that $Z$ is \emph{locally $\mathcal{P}$-contractible} if every point $z_0\in Z$ has an open neighborhood $U$ which is \emph{$\mathcal{P}$-contractible}, in the sense that there exists a homotopy $H:U\times[0,1]\to U$ between the constant function $U\ni z\mapsto z_0$ and the identity ${\rm id}_U$, which satisfies
\begin{enumerate}[(a)]
      \item $H(z_0,t)=z_0$, for all $t\in[0,1]$,
    %  \item $H(z,0)=z_0$ and $H(z,1)=z$, for all $z\in U$;
      \item for every $z\in U$, the path $p_z\coloneqq H(z,t)$ belongs to $\mathcal{P}$.
\end{enumerate}
Further, let $Z'$ also be a metric space and $p:[0,1]\to Z'$ be a path in $Z'$. We say that a continuous map $f:Z\to Z'$ has the \emph{continuation property for $p$}, if for every $b\in(0,1]$ and every path $q\in\cC([0,b), Z)$ such that $f\circ q=p\restriction_{[0,b)}$, there exists a sequence $\{t_n\}$ in $[0,b)$ convergent to $b$ and such that $\{q(t_n)\}$ converges in $Z$. Furthermore, a continuous map $f:Z\to Z'$ is called a \emph{covering map}, if every $z'\in Z'$ has an open neighborhood $U$ such that $f^{-1}(U)$ is the disjoint union of open subsets of $Z$ each of which is mapped homeomorphically into $U$ by $f$.
\end{Def}

\begin{proof}[\textit{Proof of Theorem~\ref{Thm:global-IF}}]
We proceed in several steps.

\emph{Step 1}: We first reduce the original problem to a simpler case. By the existence of $\phi,\psi$ in assumption~(iii) it follows that, for $i\in\{1,2\}$, $\pi_i(Z)$ are open sets in $\mathcal U, \mathcal V$, resp., and since $Z\subseteq \pi_1(Z)\times\pi_2(Z)$, it is no loss of generality to assume $X\times Y=\pi_1(Z)\times\pi_2(Z)$. That is, we search for an implicit function for the restriction $F:\pi_1(Z)\times\pi_2(Z)\rightarrow \mathcal{Z}$ instead of $F:X\times Y\rightarrow\mathcal{Z}$. Next, we argue that it suffices to prove the theorem for cases in which (i)--(iii) are satisfied with $\phi={\rm id}_{\mathcal{X}}$ and $\psi={\rm id}_{\mathcal{Y}}$. Note that these assumptions imply $\mathcal{U}=X=\pi_1(Z)=\mathcal{X}$ and $\mathcal{V}=Y=\pi_2(Z)=\mathcal{Y}$ since $X$ and $Y$ are open subspaces of $\mathcal{U}$ and $\mathcal{V}$, resp. Having proved this case, we can conclude the general case by considering the function $\tilde F={F}\circ(\phi^{-1},\psi^{-1})$ with $\tilde{F}:\cX \times \cY \tto \mathcal{Z}$. Next we translate the conditions (i)--(iii) on~$F$ to conditions on~$\tilde{F}$ for $\tilde{Z} \coloneqq (\phi,\psi)(Z)$.
\begin{enumerate}[(i)']
    \item We have that $\tilde{Z}$ is path-connected and closed in $\cX\times \cY$ if, and only if, $Z$ is path-connected and closed in $X\times Y$.
    \item Define
    \[
        \tilde{S}:\tilde{Z}\to \cY,\ (\tilde x,\tilde y)\mapsto \big(D(\psi^{-1})(\tilde{y})\big)^{-1}\cdot S(\phi^{-1}(\tilde{x}),\psi^{-1}(\tilde{y})).
    \]
    With the identification $(x,y)=(\phi^{-1}(\tilde{x}),\psi^{-1}(\tilde{y}))$ we obtain for all $(x,y)\in Z$ that
    \begin{align*}
      D_y F(x,y) &= D_y \big( \tilde F(\phi(x),\psi(y))\big)\\
      &= (D_{\tilde y} \tilde F)(\phi(x),\psi(y))\cdot D\psi(y)\\
      &=  D_{\tilde{y}}\tilde{F}(\tilde{x},\tilde{y}) \left(D(\psi^{-1})(\tilde{y})\right)^{-1},
    \end{align*}
    where the latter equality is a consequence of the inverse function theorem. Using this we find that
    \begin{align*}
    &\ S(x,y)D_yF(x,y)={\rm id}_\mathcal{V}
    \\\Longleftrightarrow\ \ &S\big(\phi^{-1}(\tilde{x}),\psi^{-1}(\tilde{y})\big)\cdot\left(D_{\tilde y}\tilde{F}\right)\!\!\big(\tilde x, \tilde y\big)\cdot \left(D(\psi^{-1})(\tilde{y})\right)^{-1}={\rm id}_\mathcal{V}
    \\\Longleftrightarrow\ \ & \underset{=\tilde S(\tilde x,\tilde y)}{\underbrace{\left(D(\psi^{-1})(\tilde{y})\right)^{-1} \cdot S\big(\phi^{-1}(\tilde{x}),\psi^{-1}(\tilde{y})\big)}} \cdot D_{\tilde{y}}\tilde{F}(\tilde{x},\tilde{y})={\rm id}_\mathcal{Y}.
    \end{align*}

    \item Similar to the computations above we obtain that, for $(x,y)\in Z$,
    \[
        D_{\tilde{x}}\tilde{F}(\tilde{x},\tilde{y}) = D_{{x}}{F}({x},{y})\cdot \big(D\phi({x})\big)^{-1}.
    \]
    Therefore, we have for all continuous weights $\omega:[0,\infty)\to(0,\infty)$ and all $(x,y)\in Z$ that, omitting the spaces in the subscripts of the norms,
    \begin{align*}
    &\left\| D\psi({y})\cdot{S}({x},{y})\right\|\cdot\left\|D_{{x}}{F}({x},{y})\cdot \big(D\phi({x})\big)^{-1}\right\| \le \omega(\|\psi(y)\|)
    \\ \Longleftrightarrow\ \ & \left\|\tilde{S}(\tilde{x},\tilde{y})\right\|\cdot\left\|D_{\tilde{x}}\tilde{F}(\tilde{x},\tilde{y})\right\|\leq\omega(\|\tilde{y}\|).
  \end{align*}
\end{enumerate}
Now, define the projections $\tilde{\pi}_i(p_1,p_2)=p_i$ with $i\in\{1,2\}$, $(p_1,p_2)\in\mathcal{X}\times\mathcal{Y}$. We recapitulate the situation with the following commuting diagram:
    \begin{center}
    \begin{tikzpicture}
	  \node (Z) 						  				{$X \times Y \supseteq{Z}$};
	  \node (P1)  [right of=Z, above of=Z,xshift=2ex,yshift=2ex]	{$X=\pi_1(Z)$};
	  \node (P2)  [right of=Z, below of=Z,xshift=2ex,yshift=-2ex]				{$Y=\pi_2(Z)$};
	  \node (TP1) [node distance=4cm, right of=P1]	{$\tilde{\pi}_1(\tilde{Z})\subseteq\cX$};
	  \node (TP2) [node distance=4cm, right of=P2]	{$\tilde{\pi}_2(\tilde{Z})\subseteq\cY$};
	  \node (TZ)  [node distance=6.5cm, right of=Z] 			{$\tilde Z\subseteq \mathcal{X}\times\mathcal{Y}$};
	
	  \draw[->] (Z.15)		to node [left]	{$\pi_1$} 			(P1.-60);
	  \draw[->] (Z.-15) 	to node [left]	{$\pi_2$} 			(P2.60);
	  \draw[->] (Z)  		to node [above]	{$(\phi,\psi)$} 	(TZ);
	  \draw[->] (P1) 		to node [above]	{$\phi$}			(TP1);
	  \draw[->]	(P2) 		to node [above]	{$\psi$}			(TP2);
	  \draw[->]	(TZ.165)	to node [right]	{$\tilde{\pi}_1$}	(TP1.-140);
	  \draw[->] (TZ.-165)	to node [right]	{$\tilde{\pi}_2$}	(TP2.140);
	  % \draw[->, bend right] (P1) to node [swap] {$\hat{g}$} (A);
	  % \draw[->, bend left] (P1) to node {$\hat{f}$} (B);
	  % \draw[->, dashed] (P1) to node {$k$} (P);
	\end{tikzpicture}
	\end{center}For the conclusion, consider $\tilde{g}\coloneqq\psi\circ g\circ\phi^{-1}$ together with the equality $\pi_1(Z)=\phi^{-1}(\tilde{\pi}_1(\tilde{Z}))$.

\emph{Step 2}: By Step~1, in the following we assume that $X=\pi_1(Z)$ and $Y=\pi_2(Z)$ as well as $\phi = {\rm id}_{\cX}$ and $\psi = {\rm id}_{\cY}$. We show that $\pi_1:Z\tto\pi_1(Z)$ is a local homeomorphism between connected metric spaces. Clearly,~$Z$ and~$\pi_1(Z)$ are metric spaces and since $Z$ is path-connected, $\pi_1(Z)$ is path-connected as well. To show that~$\pi_1:Z\tto\pi_1(Z)$ is a local homeomorphism, let $(a,b)\in Z$, i.e., $F(a,b)=0$. Applying the implicit function theorem, see e.g.~\cite[Thm.~10.2.1]{Dieu69}, yields open neighborhoods $U\subseteq X$ of~$a$, $V\subseteq Y$ of~$b$, and $g\in\cC^1(U, V)$ such that
\[
    \setdef{ (x,g(x))}{x\in U} = \setdef{(x,y)\in U\times V}{ S(a,b)F(x,y) = 0}.
\]
Consider the restriction $\hat \pi_1:Z\cap (U\times V)\tto\pi_1\big(Z\cap (U\times V)\big)$. Then $\hat \pi_1$ is injective since $\hat \pi_1(x_1,y_1)=\hat \pi_1(x_2,y_2)$ for some $(x_1,y_1), (x_2,y_2)\in Z\cap (U\times V)$ gives $x_1 = x_2$ and $S(a,b)F(x_1,y_1) = S(a,b)F(x_2,y_2)$, thus $y_1 = g(x_1) = g(x_2) = y_2$. Therefore, $\hat \pi_1$ is bijective and continuous. Furthermore, it is easy to see that $\hat \pi_1$ is an open map, and hence it is a homeomorphism.

\emph{Step 3}: Let
\[
    \cP \coloneqq  \cC^1\big([0,1],\pi_1(Z)\big)
\]
and observe that $\pi_1(Z)$ is $\cP$-connected and locally $\cP$-contractible since it is open. We show that~$\pi_1$ has the continuation property for every path in~$\cP$, that is, for all $q_1\in \cP$, all $b\in (0,1]$ and all $q_2\in\cC([0,b), Y)$ such that $(q_1(t),q_2(t))\in Z$ for all $t\in[0,b)$ there exists a sequence $(t_n)_{n\in\N}\subseteq [0,b)$ with $\lim_{n\to\infty} t_n = b$ such that $\big(q_2(t_n)\big)_{n\in\N}$ converges and
\[
    \lim_{n\to\infty} (q_1(t_n), q_2(t_n)) \in Z.
\]
First note that~$q_2$ is differentiable at any~$t\in[0,b)$, since there exists a local implicit function as in Step~2, so that $q_2(s) = g(q_1(s))$ for all~$s$ in a neighborhood of~$t$. Since~$g$ and~$q_1$ are differentiable we obtain $\dot q_2(t) = Dg(q_1(t)) \dot q_1(t)$. Moreover, it can be seen that the derivative is continuous at each point in~$[0,b)$. Then, using property~(iii), it can be proved by only a slight modification of the proof of~\cite[Cor.~5.3]{GutuJara07} that for any sequence $(t_n)_{n\in\N}\subseteq [0,b)$ with $\lim_{n\to\infty} t_n = b$ the sequence $\big(q_2(t_n)\big)_{n\in\N}$ is a Cauchy sequence and hence converges in~$Y=\mathcal{V}$. Since~$Z$ is closed in $X\times Y = \mathcal{U}\times\mathcal{V}$ by~(i) we thus obtain $\lim_{n\to\infty} (q_1(t_n), q_2(t_n)) \in Z$.

\emph{Step 4}: We show that $\pi_1:Z\tto\pi_1(Z)$ is a homeomorphism. By~\cite[Thm.~2.6]{GutuJara07} and Step~3 we may infer that~$\pi_1$ is a covering map. Since~$\pi_1(Z) = \phi^{-1}(\mathcal{X})$ is in particular simply connected by~(iii) it follows from~~\cite[Prop.~A.79]{Lee12} that~$\pi_1:Z\tto\pi_1(Z)$ is a homeomorphism.

\emph{Step 5}: By Step~4 we have $\left(x\mapsto (x,g(x)) = \pi_1^{-1}(x)\right)\in\cC(\pi_1(Z), Z)$ which uniquely defines the desired function $g\in\cC(\pi_1(Z), Y)$. Since~$\pi_1(Z)$ is in particular open by condition~(iii), for all $x\in\pi_1(Z)$ we have that~$g$ coincides with any solution provided by the implicit function theorem as in Step~1 in a neighborhood of~$x$. The implicit function theorem provides Fr\'echet-differentiability of the local solution, thus~$g$ is Fr\'echet-differentiable at~$x$.
\end{proof}

We like to emphasize that~$Z$ in Theorem~\ref{Thm:global-IF} may only be a subset of the zero set~$F^{-1}(0)$. This allows to exclude points~$(x,y)$ in~$F^{-1}(0)$ at which~$D_y F(x,y)$ has not left inverse, or, actually, to exclude open sets containing such points so that~$Z$ is closed (alternatively, one may restrict the sets~$X$ and~$Y$). Then a global implicit function may still exist in each connected component of~$Z$, provided the growth bound in~(iii) is satisfied.

\begin{Rem}\label{Rem:indep-diff}
An important question that arises is whether the growth bound in condition~(iii) in Theorem~\ref{Thm:global-IF} is independent of the choice of the diffeomorphisms~$\phi$ and~$\psi$. Use the notation from Theorem~\ref{Thm:global-IF}, assume that conditions~(i)--(iii) are satisfied and let $\hat\phi:\pi_1(Z)\to\hat\cX$ and $\hat\psi:\pi_2(Z)\to\hat\cY$ be diffeomorphisms for some Banach spaces~$\hat\cX, \hat\cY$. Then, omitting the subscripts indicating the spaces corresponding to the norms, we have the estimate
\begin{align*}
	& \left\| D\hat\psi({y})\cdot{S}({x},{y})\right\| \cdot\left\|D_{{x}}{F}({x},{y})\cdot \big(D\hat\phi({x})\big)^{-1}\right\|\\
    & \le \left\| D\psi({y})\cdot{S}({x},{y})\right\| \cdot \left\|D\hat\psi({y})\cdot(D\psi({y})\big)^{-1}\right\| \\
    &\quad  \cdot\left\|D_{{x}}{F}({x},{y})\cdot \big(D\phi({x})\big)^{-1}\right\| \cdot \left\|D\phi({x})\cdot(D\hat\phi({x})\big)^{-1}\right\|\\
    &\le \omega(\|\psi(y)\|_\cY) \cdot \left\|D\hat\psi({y})\cdot(D\psi({y})\big)^{-1}\right\| \cdot \left\|D\phi({x})\cdot(D\hat\phi({x})\big)^{-1}\right\|
\end{align*}
for all $(x,y)\in Z$. If the last term satisfies
\begin{multline*}
   \forall\, (x,y)\in Z:\ \omega(\|\psi(y)\|_\cY) \cdot \left\|D\hat\psi({y})\cdot(D\psi({y})\big)^{-1}\right\| \cdot \left\|D\phi({x})\cdot(D\hat\phi({x})\big)^{-1}\right\|\\ \le \hat\omega(\|\hat\psi(y)\|_{\hat \cY})
\end{multline*}
for some continuous weight~$\hat \omega$, then  the growth bound in condition~(iii) would indeed be independent of~$\phi$ and~$\psi$. However, it is still an open problem whether this is true (or a counterexample exists) and remains for future research.
\end{Rem}

%%%%
\section{Examples and discussion}\label{Sec:Ex}
%%%

In this section we discuss the assumptions in Theorem~\ref{Thm:global-IF} and provide some illustrative examples. First, we provide a practical example occurring in the modelling of electrical circuits, where the projection $\pi_1(Z)$ (on which the implicit function is defined) is a proper subset of $X$.

\begin{Ex}
Consider two diodes $\mathcal{D}_1$, $\mathcal{D}_2$ with associated currents $i_1,i_2$ and voltages $u_1,u_2$. Following \cite[Eq.\ (39.46)]{Kuepf17}, we can model their constitutive relations as
\begin{equation}\label{eq:ConstRel}i_1(t)=a_1\left(e^{\tfrac{u_1(t)}{b_1}}-1\right),\quad i_2(t)=a_2\left(e^{\tfrac{u_2(t)}{b_2}}-1\right),\quad t\in\R,\end{equation}
for some constants $a_1,a_2,b_1,b_2>0$. We may further impose the restrictions
\begin{equation}\label{eq:VoltConst}u_1(t)\in(u_1^{\min},u_1^{\max})\eqqcolon Y_1,\  u_2(t)\in(u_2^{\min},u_2^{\max})\eqqcolon Y_2,\quad t\in\R,\end{equation}
reflecting physical properties, e.g. the regions of operation of the corresponding devices that are being modelled. From \eqref{eq:ConstRel} and \eqref{eq:VoltConst} we can also derive restrictions
\begin{equation}\label{eq:CurrConst}i_1(t)\in(i_1^{\min},i_1^{\max})\eqqcolon X_1,\  i_2(t)\in(i_2^{\min},i_2^{\max})\eqqcolon X_2,\quad t\in\R.\end{equation} Such restrictions are incorporated in a natural manner in the port-Hamiltonian modelling of nonlinear electrical circuits in the context of resistive relations, see e.g.~\cite{GerHalReivdS21}. Next, we consider a parallel connection of the two diodes and add a current source with current $I$ and voltage $V$ as depicted in Figure~\ref{Fig:diodes}, i.e., we have
\begin{equation}\label{eq:Kirchhoff}
\begin{aligned}
&i_1(t)+i_2(t)=I(t),\\
&u_1(t)=u_2(t)=V(t),\quad\text{for }t\in\R.
\end{aligned}
\end{equation}
\begin{figure}
	 \centering
		\begin{circuitikz}%[every loop/.style={<-,shorten <=1pt,min distance=12mm}]
		\draw
		    (5,0) to [I,i=$I$,v_<=$V$] (1,0) to [short] (0,0) to [short] (0,3) to [short,-*] (1,3) to [short] (1,4)  to[Do,v^>=$u_{1}$,i>_=$\quad i_{1}$]  (5,4)
		    to [short,-*] (5,3) to [short] (6,3) to [short] (6,0) to [short] (5,0);
		\draw (1,3) to [short] (1,2) to[D,v^>=$u_{2}$,i>_=$\quad i_{2}$] (5,2) to [short] (5,3);
		\end{circuitikz}
		\caption{Circuit containing two diodes.}\label{Fig:diodes}
\end{figure}
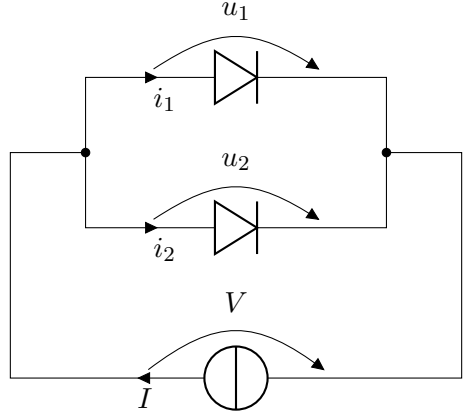

\noindent From the constitutive relations, it is clear that we can describe $i_1,i_2$ and hence $I$ as a function of $V$. However, the converse is not evident, i.e., how $V$ is given in terms of the current $I$. Invoking \eqref{eq:ConstRel} and \eqref{eq:Kirchhoff},  $I$ and $V$ satisfy the relation
\[I=a_1\left(e^{\tfrac{V}{b_1}}-1\right)+a_2\left(e^{\tfrac{V}{b_2}}-1\right)\eqqcolon f(V).\]
Further by \eqref{eq:VoltConst} and \eqref{eq:Kirchhoff}, $V$ has to satisfy $V\in Y_1\cap Y_2\eqqcolon Y$, whereas $I$ has to satisfy $I\in X_1 + X_2 = (i_1^{\min}+i_2^{\min}, i_1^{\max}+i_2^{\max})\eqqcolon X$ by \eqref{eq:CurrConst} and \eqref{eq:Kirchhoff}. Note that both $X$ and $Y$ are open intervals and we exclude the trivial case that $Y$ is empty. %Furthermore, we assume that one can only apply a certain amount of current, i.e., $I\in [I_*,I^*]$ for some $I_*,I^*\in\R$.
Recapitulating, with $F:X\times Y\rightarrow\R, (I,V)\mapsto I-f(V)$ and
\[Z\coloneqq \setdef{(I,V)\in X\times Y}{F(I,V)=0},\]
we seek the existence of a function $g\in\cC(\pi_1(Z),Y)$ such that
\[\setdef{(I,g(I))}{I\in\pi_1(Z)}=Z.\]
The existence of such a function is obviously equivalent to the invertibility of $f$ on $\pi_2(Z)\subseteq Y$, which holds true since its derivative is positive. Nevertheless, we check the assumptions (i)-(iii) of Theorem \ref{Thm:global-IF} in order to illustrate it. Since $Z$ is the zero set of a continuous function, it is relatively closed in $X\times Y$, i.e., (i) holds. For (ii), note that
\[D_VF(I,V)=-f'(V)=-\tfrac{a_1}{b_1}\exp\left(\tfrac{V}{b_1}\right)-\tfrac{a_2}{b_2}\exp\left(\tfrac{V}{b_2}\right)<0.\]
It remains to find diffeomorphisms $\phi:\pi_1(Z)\tto\R$, $\psi:\pi_2(Z)\tto\R$ and a continuous weight $\omega:[0,\infty)\to(0,\infty)$ such that the growth bound in~(iii) is satisfied. Choose any diffeomorphism $\phi:\pi_1(Z)\rightarrow\R$, which exists since $\pi_1(Z) = X$ is an open interval. Define $\psi\coloneqq\phi\circ f$, which is a diffeomorphism since $f$ is invertible on $\pi_2(Z)$. Then $\psi'=(\phi'\circ f)\cdot f'$. Further, let $\omega(t)=t+1$ for $t\in[0,\infty)$ and note that $S(I,V) = \big(D_V F(I,V)\big)^{-1} = \big(-f'(V)\big)^{-1}$ for $(I,V)\in Z$. Recalling that $I=f(V)$ for all $(I,V)\in Z$ we find
\begin{align*}
	&\left\| D\psi({V})\cdot{S}({I},{V})\right\| \cdot\left\|D_{{I}}{F}({I},{V})\cdot \big(D\phi({I})\big)^{-1}\right\|\\
	& = \left| \phi'(I)\cdot f'(V)\cdot \big(-f'(V)\big)^{-1}\right| \cdot\left|\phi'({I})^{-1}\right| = \frac{|\phi'({I})|}{|\phi'({I})|} = 1\\
	& \le|\psi(V)|+1 = \omega(\|\psi(V)\|),
\end{align*}
for all $(I,V)\in Z$, proving (iii).
\end{Ex}

We like to highlight that none of the assumptions (i)--(iii) in Theorem~\ref{Thm:global-IF} can be omitted in general. It is clear that connectedness of~$Z$ in~(i) and local solvability guaranteed by~(ii) are indispensable. Counterexamples in finite dimension are constructed for~(iii) in the following examples. Condition~(iii) basically consists of two parts. The first one is to check whether~$\pi_i(Z)$, $i=1,2$, are diffeomorphic to some Banach spaces. The second part is the growth bound involving the diffeomorphisms, the partial derivative~$D_x F$ and the left inverse of~$D_y F$.

First, we like to discuss why we chose the projections of~$Z$ as the domains of the diffeomorphisms in Theorem~\ref{Thm:global-IF}, whereas intuitively one could consider the open sets~$X$ and~$Y$ as the domains.

\begin{Rem}
In a possible different formulation of Theorem~\ref{Thm:global-IF} one could choose diffeomorphisms $\tilde\phi:X\tto\cX$ and $\tilde\psi:Y\tto\cY$ and then consider, \textit{mutatis mutandis}, the corresponding growth bound in condition~(iii). This would relax the assumptions on the projections~$\pi_i(Z)$, which would then not necessarily need to be open and simply connected. However, for the proof technique to be feasible we need to additionally require that~$\pi_1(Z)$ is simply connected. Indeed, the proof is analogous, but the modified theorem does not cover basic examples.

For $F: \R\times (-1,1)\to\R,\ (x,y)\mapsto x - y$ and $Z:= F^{-1}(0)$ it is easy to check that conditions~(i) and~(ii) are satisfied. The growth bound in condition~(iii) reads
\begin{align*}
                & \forall\, (x,y)\in Z:\ |\tilde\psi'(y)|\cdot|\tilde\phi'(x)^{-1}| \le \omega(|\tilde\psi(y)|) \\
 \iff\ \ & \forall\, y\in(-1,1):\ |\tilde\phi'(y)|\geq\frac{|\tilde\psi'(y)|}{\omega(|\tilde\psi(y)|)}
\end{align*}
for some continuous weight~$\omega$. Note that $\tilde\phi'((-1,1))$ is bounded and $\tilde\phi(y)\neq 0$ for all $y\in\R$, hence
\[
    \int_{-1}^1|\tilde\phi'(y)| {\rm d}y = \left|\int_{-1}^1 \tilde\phi'(y)  {\rm d}y\right| < \infty.
\]
Then the change of variables theorem with the substitution $t=\tilde\psi(y)$ together with the inverse function theorem yields that
\begin{align*}
     \infty &> \int_{-1}^1|\tilde\phi'(y)| {\rm d}y \ge \int_{-1}^1 \frac{|\tilde\psi'(y)|}{\omega(|\tilde\psi(y)|)} {\rm d}y = \int_{-\infty}^\infty \frac{\left|\tilde\psi'\big(\tilde\psi^{-1}(t)\big)\right|}{\omega(|t|)} \left|(\tilde \psi^{-1})'(t)\right| {\rm d}t \\ &= 2 \int_0^\infty \frac{1}{\omega(t)} {\rm d}t = \infty,
\end{align*}
a contradiction.

Nevertheless, a global implicit function obviously exists. The assumptions of Theorem~\ref{Thm:global-IF} are satisfied since $\pi_1(Z) = \pi_2(Z) = (-1,1)$ and we may choose $\phi = \psi$, with which the growth bound holds true.
\end{Rem}

%\begin{Ex}
%  Consider
%  \[
%    F: \R\times (-1,1)\to\R,\ (x,y)\mapsto x - y,
%  \]
%  for which it is easy to check that conditions~(i) and~(ii) in Theorem~\ref{Thm:global-IF} are satisfied. For condition~(iii) we observe that $\pi_1(Z) = \pi_2(Z) = (-1,1)$ and we choose $\phi = \psi$ with
%  \[
%    \phi:(-1,1) \to \R,\ x\mapsto \tan\left(\tfrac{\pi}{2} x\right).
%  \]
%  Since $(x,y)\in Z$ implies $x=y$ we find that $|D\psi(y)| / |D\phi(x)| = 1$ and hence the growth bound in condition~(iii) is satisfied for~$\omega(t) = 1 + t$ for instance.
%\end{Ex}

We continue by presenting an example where in assumption~(iii) it is not possible to find suitable diffeomorphisms and, at the same time, a global implicit function does not exist.

\begin{Ex}
Consider
\[
    F:\R^2\times\R\to \R^2,\ (x_1,x_2,y)\mapsto \begin{pmatrix} x_1 -\cos y\\ x_2 - \sin y\end{pmatrix},\quad Z:=  F^{-1}(0).
\]
Then assumptions~(i) and (ii) in Theorem~\ref{Thm:global-IF} are satisfied. Since $\pi_1(Z)=S^1$, the unit circle in~$\R^2$, there is no Banach space~$\cX$ such that~$\pi_1(Z)$ is diffeomorphic to~$\cX$. Indeed, no global implicit function can exist, since $y\mapsto (\cos y, \sin y)$ is not injective on $\R$.
% We have $\phi = \psi={\rm id}$ and condition~(iv) translates to
% \begin{align*}
% &\left\|\left(\left(\frac{\partial F}{\partial y}(x,y)\right)^\top\frac{\partial F}{\partial y}(x,y)\right)^{-1}\left(\frac{\partial F}{\partial y}(x,y)\right)^\top\right\|\cdot\left\| D_xF(x,y)\right\|
% \\&\leq \omega(|y|)
% \\\Longleftrightarrow\quad&\left\| \begin{bmatrix}-\sin y&\cos y\end{bmatrix}\right\|\leq \omega (|y|).
% \end{align*}
% Note that we may choose any norm in~(iv) by the equivalence of all norms in $\R^n$. Furthermore, if $\omega$ is a continuous weight, then $c_1\omega(c_2\cdot)$ is also a continuous weight for every $c_1,c_2\in\R$. Hence, choosing $\|\cdot\|$ as the Euclidean norm and the weight $\omega(t)=1+t$, we obtain the inequality \[1\leq 1+|y|\] for (iv), which is true for every $y\in\R$.
% Now, note that $\pi(Z) = \setdef{x\in\R^2}{\|x\|=1}$, which is not simply connected, i.e.,~(iii) is not satisfied.
\end{Ex}

% Indeed, the equation $F(x,y)=0$ is just the real-valued version of the complex-valued equation
% \[
%     x - \exp(y) = 0,\quad x\in \C\setminus\{0\},\ y\in \C,
% \]
% the solution of which is the complex logarithm, which cannot be represented as a continuous function $g:\C\setminus\{0\}\tto\C$, i.e., $F(x,y)=0$ does not have a unique maximal solution.\\

In the next example the growth bound in condition~(iii) is not satisfied for any suitable choice of diffeomorphisms and, at the same time, a global implicit function does not exist.

\begin{Ex}
We choose $F:X\times Y\tto\R^2=\mathcal{Z}$ as a function of the type $F(x,y)=x-\tilde{F}(y)$ and $Z:= F^{-1}(0)$. This means that the existence of a global implicit function is equivalent to $\tilde{F}$ being injective. We further set $X\times Y=\R^2\times(0,1)^2$ and construct~$\tilde{F}$ by successively defining its restrictions $\tilde{F}\restriction_{(0,\delta)\times(0,1)} = \tilde{F}_1$ and $\tilde{F}\restriction_{(\delta,1)\times(0,1)} = \tilde{F}_2$ for some $0<\delta\leq\tfrac{1}{2}$. Choose $\varepsilon, \alpha>0$, and consider the non-injective function
\[
    \tilde{F}_1: (0,\delta)\times(0,1) \tto\R^2,\ (y_1,y_2)\mapsto \begin{pmatrix}(\alpha+y_2)\sin\left(\tfrac{2\pi}{\delta}(1+\varepsilon)y_1\right)\\ (\alpha+y_2)\cos\left(\tfrac{2\pi}{\delta}(1+\varepsilon)y_1\right)\end{pmatrix}.
\]
Observe that $\im \tilde{F}_1=B_{\alpha+1}(0)\setminus\overline{B_{\alpha}(0)}$, where $B_\alpha(z)$ denotes the open ball with radius~$\alpha$ around $z\in\R^2$, i.e., the image of~$\tilde{F}_1$ is an annulus. Next, define $\tilde{F}_2$ similarly to $\tilde{F}_1$ using elementary functions such that the hole of the annulus is filled as displayed in Fig.~\ref{FIG1}, i.e., $\overline{B_\alpha(0)}\subset\im\tilde{F}_2\subset B_{\alpha+1}(0)$. Note that $\tilde{F}_2$ can be chosen such that the resulting composition~$\tilde{F}$ is differentiable everywhere.  Overall, we have constructed a \emph{non-injective} function~$\tilde{F}$.

\def\svgwidth{\columnwidth}
\begingroup%
  \makeatletter%
  \providecommand\color[2][]{%
    \errmessage{(Inkscape) Color is used for the text in Inkscape, but the package 'color.sty' is not loaded}%
    \renewcommand\color[2][]{}%
  }%
  \providecommand\transparent[1]{%
    \errmessage{(Inkscape) Transparency is used (non-zero) for the text in Inkscape, but the package 'transparent.sty' is not loaded}%
    \renewcommand\transparent[1]{}%
  }%
  \providecommand\rotatebox[2]{#2}%
  \newcommand*\fsize{\dimexpr\f@size pt\relax}%
  \newcommand*\lineheight[1]{\fontsize{\fsize}{#1\fsize}\selectfont}%
  \ifx\svgwidth\undefined%
    \setlength{\unitlength}{287.38166869bp}%
    \ifx\svgscale\undefined%
      \relax%
    \else%
      \setlength{\unitlength}{\unitlength * \real{\svgscale}}%
    \fi%
  \else%
    \setlength{\unitlength}{\svgwidth}%
  \fi%
  \global\let\svgwidth\undefined%
  \global\let\svgscale\undefined%
  \makeatother%
  \begin{figure}
  \begin{picture}(1,0.43719072)%
    \lineheight{1}%
    \setlength\tabcolsep{0pt}%
    \put(0,0){\includegraphics[width=\unitlength,page=1]{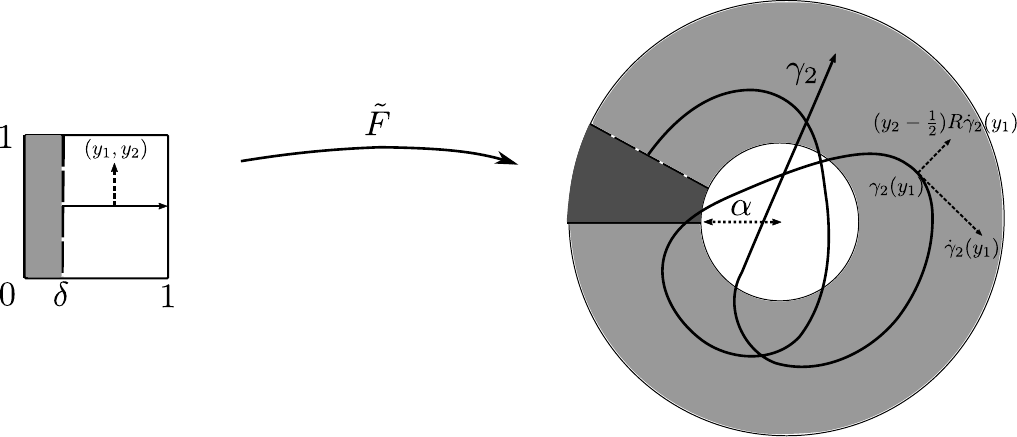}}%
    %\put(0.35390194,0.30828789){\color[rgb]{0,0,0}\makebox(0,0)[lt]{\lineheight{1.25}\smash{\begin{tabular}[t]{l}$\tilde{F}$\end{tabular}}}}%
%    \put(0.01896377,0.12877993){\color[rgb]{0,0,0}\makebox(0,0)[t]{\lineheight{1.25}\smash{\begin{tabular}[t]{c}$0$\end{tabular}}}}%
%    \put(0.0423607,0.12891283){\color[rgb]{0,0,0}\makebox(0,0)[lt]{\lineheight{1.25}\smash{\begin{tabular}[t]{l}$\delta$\end{tabular}}}}%
%    \put(0.15972444,0.12901307){\color[rgb]{0,0,0}\makebox(0,0)[t]{\lineheight{1.25}\smash{\begin{tabular}[t]{c}$1$\end{tabular}}}}%
%    \put(0.01836013,0.31551138){\color[rgb]{0,0,0}\makebox(0,0)[t]{\lineheight{1.25}\smash{\begin{tabular}[t]{c}$1$\end{tabular}}}}%
%    \put(0,0){\includegraphics[width=\unitlength,page=2]{FIG1_new.pdf}}%
%    \put(0.73325068,0.22207236){\color[rgb]{0,0,0}\makebox(0,0)[t]{\lineheight{1.25}\smash{\begin{tabular}[t]{c}$\alpha$\end{tabular}}}}%
  \end{picture}%
  \caption{Illustration of the construction of $\tilde{F}$.}
  \label{FIG1}
  \end{figure}
\endgroup%

%\ \\
%\textbf{Im Bild noch Pfeil mit $\gamma_2$ auf die Kurve einfügen}\\

Observe that $\im \tilde F = \pi_1(Z) = B_{\alpha+1}(0)$ and $\pi_2(Z) = (0,1)^2$. Then the three conditions on $F$ translate to $\tilde{F}$ as follows:
\begin{enumerate}[(i')]
    \item the graph of $\tilde{F}$ is connected;
    \item $\forall\, y\in (0,1)^2:\ \rk D\tilde{F}(y) = 2$\,;
    \item there exist diffeomorphisms $\phi:B_{\alpha+1}(0)\rightarrow \mathcal{X},\psi:(0,1)^2\rightarrow\mathcal{Y}$ for some Banach spaces $(\cX,\|\cdot\|_\cX), (\cY,\|\cdot\|_\cY)$, and a continuous weight $\omega:[0,\infty)\to(0,\infty)$ such that for all $y\in (0,1)^2$ we have
\[
        \left\| D\psi(y) \cdot D\tilde{F}_1(y)^{-1}\right\|_{\mathcal L(\mathcal Z,\mathcal Y)} \cdot \left\|\big(D\phi(\tilde F_1(y))\big)^{-1}\right\|_{\mathcal L(\mathcal X,\mathcal Y)} \le \omega(\|\psi(y)\|_{\mathcal Y}).
\]
\end{enumerate}
Note that~(i') is guaranteed by our choice of $Y=(0,1)^2$ and the continuity of~$\tilde{F}$, which holds by construction. For~(ii') note that
\[
    D\tilde{F}_1(y_1,y_2)= \begin{bmatrix}(\alpha+y_2)\tfrac{2\pi}{\delta}(1+\varepsilon)\cos\left(\tfrac{2\pi}{\delta}(1+\varepsilon)y_1\right)& \sin\left(\tfrac{2\pi}{\delta}(1+\varepsilon)y_1\right)\\ -(\alpha+y_2)\tfrac{2\pi}{\delta}(1+\varepsilon)\sin\left(\tfrac{2\pi}{\delta}(1+\varepsilon)y_1\right)& \cos\left(\tfrac{2\pi}{\delta}(1+\varepsilon)y_1\right)\end{bmatrix},
\]
and $\det\left(\tilde{F}_1(y_1,y_2)\right)=(\alpha+y_2)\tfrac{2\pi}{\delta}(1+\varepsilon)\neq0$. Further, the use of the constant $\alpha$ (large enough) guarantees that $\rk D\tilde{F}(y)=2$ when filling $\overline{B_\alpha(0)}$ as displayed in Fig~\ref{FIG1}. Hence, (ii') is satisfied. Next, we show that (iii') is not satisfied, although, obviously, both~$\pi_1(Z)= B_{\alpha+1}(0)$ and~$\pi_2(Z)= (0,1)^2$ are diffeomorphic to some Banach spaces~$\cX$ and~$\mathcal{Y}$. Without loss of generality, we may assume that $\cX=\mathcal{Y}=\R^2$.

% First, choose the diffeomorphisms
% \[
% 	\psi:(0,1)^2\tto\R^2,\ (y_1,y_2)\mapsto\left(\psi_1(y_1),\psi_1(y_2)\right)
% \]
% with $\psi_1(x)=\tfrac{2x-1}{1-(2x-1)^2}$ and
% \[
% 	\phi:B_{\alpha+1}(0)\tto\R^2,\ x\mapsto \tfrac{x}{(\alpha+1)^2-\|x\|_2^2}.
% \]

We show that condition (iii') is not satisfied for any diffeomorphisms $\phi:B_{\alpha+1}(0)\tto\R^2$ and $\psi:(0,1)^2\tto\R^2$ by considering two cases. Let $(\hat y_1,\hat y_2) := \psi^{-1}(0,0)\in(0,1)^2$.

\emph{Case 1}: Assume that $\hat y_1 \le \delta$. We show that the growth bound fails for~$\tilde{F}_1$. Seeking a contradiction, assume that we have
\[
        \left\| D\psi(y) \cdot D\tilde{F}_1(y)^{-1}\right\|_{\mathcal L(\mathcal Z,\mathcal Y)} \cdot \left\|\big(D\phi(\tilde F_1(y))\big)^{-1}\right\|_{\mathcal L(\mathcal X,\mathcal Y)} \le \omega(\|\psi(y)\|_{\mathcal Y})
\]
for all $y\in (0,\delta)\times(0,1)$ and some weight $\omega$. Although we did not specify the norms on $\mathcal{X,Y,Z}$, using the weight property and the equivalence of all norms on $\R^{n^2},\R^n$, respectively, guarantees the existence of positive constants $c_1,c_2$ such that
\begin{align*}
&c_1\left\| D\psi(y) \cdot D\tilde{F}_1(y)^{-1}\right\|_F \cdot \left\|\big(D\phi(\tilde F_1(y))\big)^{-1}\right\|_{\mathcal L(\mathcal X,\mathcal Y)}
\\&\leq \left\| D\psi(y) \cdot D\tilde{F}_1(y)^{-1}\right\|_{\mathcal L(\mathcal Z,\mathcal Y)} \cdot \left\|\big(D\phi(\tilde F_1(y))\big)^{-1}\right\|_{\mathcal L(\mathcal X,\mathcal Y)}
\\&\leq \omega(\|\psi(y)\|_{\mathcal Y}) \leq \omega(c_2\|\psi(y)\|_{2}),
\end{align*}
where $\|\cdot\|_F$ is the Frobenius norm. Observing that $\tilde \omega(\cdot):= c_1^{-1}\omega(c_2\,\cdot)$ again defines a weight, we obtain
\begin{equation*}%\label{eq:Normen}
\left\| D\psi(y) \cdot D\tilde{F}_1(y)^{-1}\right\|_F \cdot \left\|\big(D\phi(\tilde F_1(y))\big)^{-1}\right\|_{\mathcal L(\mathcal X,\mathcal Y)}\leq \tilde\omega(\|\psi(y)\|_{2})
\end{equation*}
for all $y\in (0,\delta)\times(0,1)$. In order to simplify the computations, choose $\delta=\tfrac{1}{2}$ and $\alpha=\varepsilon=1$. Since $\tilde{F}_1\left((0,\tfrac{1}{2})\times\{\hat y_2\}\right)= (1+\hat y_2) S^1$ is a compact subset of $B_{2}(0)$ and $B_{2}(0)\ni z\mapsto \left\|D\phi(z)\right\|$ is a continuous mapping we have
\[
     \exists\, \gamma>0\ \forall\, y\in(0,\tfrac{1}{2})\times\{\hat y_2\}:\ \left\|D\phi(\tilde{F}_1(y))\right\|_{\mathcal L(\mathcal X,\mathcal Y)}\le \gamma.
\]
This gives $\left\|\big(D\phi(\tilde F_1(y))\big)^{-1}\right\|_{\mathcal L(\mathcal X,\mathcal Y)}\geq\gamma^{-1}$ for all $y\in(0,\tfrac{1}{2})\times\{\hat y_2\}$. Accordingly, we may calculate that for all $y_1\in(0,\tfrac12)$ we have
\[
D\tilde{F}_1(y_1,\hat y_2)^{-1}=\begin{bmatrix}\tfrac{1}{8\pi(1+\hat y_2)} \cos\left(8\pi y_1\right)& -\tfrac{1}{8\pi(1+\hat y_2)}\sin\left(8\pi y_1\right) \\ \sin\left(8\pi y_1\right)&\cos\left(8\pi y_1\right)\end{bmatrix}
\]
% \begin{align*}
% &D\phi(x_1,x_2)^{-1}=\frac{\big((\alpha+1)^2-\|x\|_2^2\big)^2}{(\alpha+1)^4-(x_1^4+x_2^4)}\begin{bmatrix}(\alpha+1)^2-x_1^2+x_2^2&-2x_1x_2\\-2x_1x_2&(\alpha+1)^2-x_2^2+x_1^2\end{bmatrix},
% \\&D\psi(y_1,y_2)=\begin{bmatrix}\psi_1'(y_1)&0\\0&\psi_1'(y_2)\end{bmatrix}=\begin{bmatrix}\tfrac{2}{(1-(2y_1-1)^2)^2}&0\\0&\tfrac{2}{(1-(2y_2-1)^2)^2}\end{bmatrix},
% \\&D\tilde{F}_1(y_1,y_2)^{-1}=\begin{bmatrix}d_1 \cos\left(\tfrac{2\pi}{\delta}(1+\varepsilon)y_1\right)& -d_1\sin\left(\tfrac{2\pi}{\delta}(1+\varepsilon)y_1\right) \\ \sin\left(\tfrac{2\pi}{\delta}(1+\varepsilon)y_1\right)&\cos\left(\tfrac{2\pi}{\delta}(1+\varepsilon)y_1\right)\end{bmatrix},
% \end{align*}
% where $d_1 := \tfrac{\delta}{2\pi(\alpha+y_2)(1+\varepsilon)}$.
and hence
\begin{align*}
\left\| D\psi(y_1,\hat y_2) \cdot D\tilde{F}_1(y_1,\hat y_2)^{-1}\right\|_F &=  \sqrt{\tfrac{1}{64\pi^2 (1+\hat y_2)^2}\left(\tfrac{\partial \psi_1}{\partial y_1}^2+\tfrac{\partial \psi_2}{\partial y_1}^2\right)+\tfrac{\partial \psi_1}{\partial y_2}^2+\tfrac{\partial \psi_2}{\partial y_2}^2}\\
&\ge \tfrac{1}{8\pi(1+\hat y_2)}\left\|\tfrac{\partial \psi}{\partial y_1}(y_1,\hat y_2)\right\|_2.
\end{align*}
Note that for all $y_1\in (0,\hat y_1)$ we have that $\|\psi(y_1,\hat y_2)\|_2 > 0$ and, because $\lim_{y_1\to 0} \|\psi(y_1,\hat y_2)\|_2 = \infty$, the set
\[
    \cS := \setdef{y_1\in(0,\hat y_1)}{\tfrac{\partial}{\partial y_1}\|\psi(y_1,\hat y_2)\|^2_2 < 0}
\]
has compact complement $(0,\hat y_1)\setminus\cS$. Furthermore, for all $y_1\in (0,\hat y_1)$ we have that
\begin{align*}
    \tfrac12 \left|\tfrac{\partial}{\partial y_1}\|\psi(y_1,\hat y_2)\|^2_2\right| &= \left|\psi(y_1,\hat y_2)^\top \tfrac{\partial \psi}{\partial y_1}(y_1,\hat y_2)\right|\\
    &\le \|\psi(y_1,\hat y_2)\|_2 \|\tfrac{\partial \psi}{\partial y_1}(y_1,\hat y_2)\|_2\\
    &\le 8\pi(1+\hat y_2)  \|\psi(y_1,\hat y_2)\|_2 \left\| D\psi(y_1,\hat y_2) \cdot D\tilde{F}_1(y_1,\hat y_2)^{-1}\right\|_F \\
    &\le 8\pi (1+\hat y_2) \gamma  \|\psi(y_1,\hat y_2)\|_2\, \tilde\omega(\|\psi(y_1,\hat y_2)\|_{2}).
\end{align*}
With
\[
    \xi := \int_{(0,\hat y_1)\setminus\cS} \tfrac{\tfrac{\partial}{\partial y_1}\left\|\psi\left(y_1,\hat y_2\right)\right\|^2_2}{ \|\psi(y_1,\hat y_2)\|_2\, \tilde\omega\left(\left\|\psi\left(y_1,\hat y_2\right)\right\|_{2}\right)}\, \mathrm{d}y_1 < \infty
\]
and the substitutions $t=\|\psi(y_1,\hat y_2)\|_2^2$ and $u = \sqrt{t}$ we may then derive
\begin{align*}
16\pi\gamma\hat y_1 &\geq \int_0^{\hat y_1} \tfrac{\left|\tfrac{\partial}{\partial y_1}\left\|\psi\left(y_1,\hat y_2\right)\right\|^2_2\right|}{ \|\psi(y_1,\hat y_2)\|_2\, \tilde\omega\left(\left\|\psi\left(y_1,\hat y_2\right)\right\|_{2}\right)}\, \mathrm{d}y_1
\geq \int_\cS \tfrac{\left|\tfrac{\partial}{\partial y_1}\left\|\psi\left(y_1,\hat y_2\right)\right\|^2_2\right|}{ \|\psi(y_1,\hat y_2)\|_2\, \tilde\omega\left(\left\|\psi\left(y_1,\hat y_2\right)\right\|_{2}\right)}\, \mathrm{d}y_1 \\
&=\int_\cS \tfrac{-\tfrac{\partial}{\partial y_1}\left\|\psi\left(y_1,\hat y_2\right)\right\|^2_2}{ \|\psi(y_1,\hat y_2)\|_2\, \tilde\omega\left(\left\|\psi\left(y_1,\hat y_2\right)\right\|_{2}\right)}\, \mathrm{d}y_1
- \xi + \xi\\
&=\int_0^{\hat y_1} \tfrac{-\tfrac{\partial}{\partial y_1}\left\|\psi\left(y_1,\hat y_2\right)\right\|^2_2}{ \|\psi(y_1,\hat y_2)\|_2\, \tilde\omega\left(\left\|\psi\left(y_1,\hat y_2\right)\right\|_{2}\right)}\, \mathrm{d}y_1  + \xi\\
&= \int_{\left\|\psi\left(\hat y_1,\hat y_2\right)\right\|^2_2}^\infty \tfrac{1}{\sqrt{t}\tilde\omega(\sqrt{t})}\, \mathrm{d}t +\xi = 2 \int_0^\infty \tfrac{1}{\tilde\omega(u)}\, \mathrm{d}u + \xi=\infty,
\end{align*}
a contradiction.

% \begin{align*}
%  &\left\| D\psi(y_1,y_2) \cdot D\tilde{F}_1(y_1,y_2)^{-1}\right\|_F =  \big( d_1^2 \psi_1'(y_1)^2 + 4\big)^2,\\
%  & \left\|\big(D\phi(\tilde F_1(y))\big)^{-1}\right\|_{\mathrm{vec}_1} \ge \frac{2\big( (\alpha+1)^2 - (\alpha+\tfrac12)^2\big)^2}{(\alpha+1)^4} =: d_2,\\
%  & \|\psi(y_1,y_2)\|_{\infty} = \max\big\{ |\psi_1(y_1)|, |\psi_1(\tfrac12)| = 0\big\} = |\psi_1(y_1)|,
% \end{align*}
% \[
%     d_1 d_2 |\psi_1'(y_1)| \le d_2 \big( d_1^2 \psi_1'(y_1)^2 + 4\big)^2 \le \tilde\omega(|\psi_1(y_1)|).
% \]

% Invoking $|\psi_1(y_1)| = -\psi_1(y_1) > 0$ and $\psi_1'(y_1)>0$ we obtain
% \[
% \delta \geq \int_0^\delta \tfrac{d_1 d_2 |\psi'(y_1)|}{\tilde\omega(|\psi(y_1)|)}\, \mathrm{d}y_1=-\int_0^\delta\tfrac{-d_1 d_2\psi_1'(y_1)}{\tilde\omega(-\psi_1(y_1))}\, \mathrm{d}y_1\stackrel{t=\psi_1(y_1)}{=} \int_{-\psi_1(\delta)}^\infty \tfrac{d_1 d_2}{\tilde\omega(t)}\, \mathrm{d}t=\infty,
% \]
% a contradiction. We also stress that, in view of Remark~\ref{Rem:indep-diff}, for any other choice of diffeomorphisms the growth bound in Theorem~\ref{Thm:global-IF}\,(iii) still cannot be satisfied for any continuous weight.
%This can be seen by repeating the reasoning of \emph{Step~1} of the proof of Theorem~\ref{Thm:global-IF} with the change of variables $(\hat\phi\circ\phi^{-1},\hat\psi\circ\psi^{-1})$. This also shows in the context of Theorem \ref{Thm:global-IF} that condition~(iii) may only be checked for a single pair of diffeomorphism $(\phi,\psi)$, if found.
% \input{FIG1.pdf_tex}\label{FIG 1}

\emph{Case 2}: Assume that $\hat y_1 > \delta$. We show that the growth bound fails for~$\tilde{F}_2$, which is similar to Case~1. To this end, we render the definition of~$\tilde F_2$ more precisely. First define the curve
\[
    \gamma_1:(0,\delta)\to\R^2,\ t\mapsto \begin{pmatrix}(\alpha+\tfrac12)\sin\left(\tfrac{2\pi}{\delta}(1+\varepsilon)t\right)\\ (\alpha+\tfrac12)\cos\left(\tfrac{2\pi}{\delta}(1+\varepsilon)t\right)\end{pmatrix}
\]
and the (rotation) matrix $R:= \begin{bmatrix} 0&-1\\ 1& 0\end{bmatrix}$. Then~$\tilde F_1$ can alternatively be written as
\[
    \tilde F_1(y_1,y_2) = \gamma_1(y_1) + (y_2 - \tfrac12) R \dot \gamma_1(y_1)
\]
for all $(y_1,y_2)\in(0,\delta)\times(0,1)$. In view of this, we may choose a curve~$\gamma_2:(\delta,1)\to\R^2$ as depicted in Fig.~\ref{FIG1} such that
\[
    \tilde F_2(y_1,y_2) = \gamma_2(y_1) + (y_2 - \tfrac12) R \dot \gamma_2(y_1)
\]
for all $(y_1,y_2)\in(\delta,1)\times(0,1)$ and, as mentioned before, $\overline{B_\alpha(0)} \subset \im \tilde F_2 \subset B_{\alpha+1}(0)$ and $D\tilde F_2(y_1,y_2)$ is invertible for all $(y_1,y_2)\in(\delta,1)\times(0,1)$. Omitting the details, the same arguments as in Case~1 may now be applied to arrive at a contradiction. In particular, fixing~$y_2 = \hat y_2$ leads to a curve $\tilde\gamma_2(y_1) = \gamma_2(y_1) + (\hat y_2 - \tfrac12) R \dot \gamma_2(y_1)$ along which the growth bound is violated.
\end{Ex}

\begin{Rems}\label{Rem:relax}\
Finally, we like to point out that, while Theorem~\ref{Thm:global-IF} is already quite general, still it does not cover all relevant cases. Consider
\[
    F:\R^2\times \R^2\to \R^3,\ (x_1,x_2,y_1,y_2)\mapsto \begin{pmatrix} x_1 - y_1\\ x_2 - y_2\\ x_1^2 + x_2^2 - 1\end{pmatrix},
\]
then
\[
    Z := F^{-1}(0) = \setdef{(x_1,x_2,y_1,y_2)\in\R^4}{ x_1^2 + x_2^2 = 1,\ y_1=x_1,\ y_2=x_2}
\]
and $\pi_1(Z),\pi_2(Z)$ are both the unit circle in~$\R^2$, i.e., closed subsets which are not simply connected, for which it is not possible to satisfy condition~(iii). However, a global implicit function obviously exists. Further research is necessary to cover examples of this type.
\end{Rems}

% \textbf{Problem 2}: Even if we do this, we get another problem for $X=Y=\R^2$ and $F(x_1,x_2,y_1,y_2) = (x_1-y_1, x_2-y_2, x_1)$ for which
% \[
%     Z = \setdef{(x,y)\in\R^4}{x_1 = y_1 = 0,\ x_2 = y_2},
% \]
% and then a diffeomorphism $\phi:\pi_1(Z)\to \R^2$ obviously does not exist. However, all other conditions are satisfied and a global implicit function exists.
%
%%%%%%%%%%%%%%%%%%%%%%%%%%%%%%%%%%%%%%%%%%%%%%%%%%%%%%%%%%%%%%%%%%%%%%%%%%%%%%%%%%%%%%%%%%%%%%%%%%%%%%%%%%%%%
%\section{Conclusion}\label{Sec:Concl}
%%%%%%%%%%%%%%%%%%%%%%%%%%%%%%%%%%%%%%%%%%%%%%%%%%%%%%%%%%%%%%%%%%%%%%%%%%%%%%%%%%%%%%%%%%%%%%%%%%%%%%%%%%%%%
%

\section*{Acknowledgment}

We thank Henrik Schumacher (RWTH Aachen) for several helpful discussions.

\bibliographystyle{plain}
% \bibliography{MST}
%\bibliography{Haller-database}

\end{document}